\documentclass[11pt]{amsart}

\usepackage{geometry}
\usepackage{soul}

\usepackage{amsfonts, amssymb, amscd}
\numberwithin{equation}{section}

\usepackage{bm}
\usepackage{verbatim}
\usepackage{mathrsfs}
\usepackage{graphicx}
\usepackage{tikz-cd}
\usepackage{subcaption}
\usepackage{listings}
\usepackage{subfiles}
\usepackage[toc,page]{appendix}
\usepackage{mathtools}
\usepackage{comment}
\usepackage{enumerate}
\usepackage{enumitem}
\usepackage[all]{xy}

\usepackage{graphicx}
\graphicspath{{images/}}

\usepackage{appendix}
\usepackage{hyperref}
\hypersetup{
    colorlinks=true,
    citecolor=red,
    linkcolor=blue,
    filecolor=magenta,      
    urlcolor=red,
}
\newcommand{\Qq}{\mathbb{Q}}

\newcommand{\Rr}{\mathbb{R}}

\newcommand{\vol}{\operatorname{vol}}

\newcommand{\Exc}{\operatorname{Exc}}

\newcommand{\rk}{\operatorname{rank}}

\newcommand{\ninv}{{\operatorname{ninv}}}

\newcommand{\Var}{{\operatorname{Var}}}

\newcommand{\fol}{\operatorname{fol}}

\newcommand{\pet}{\operatorname{pet}}
\newcommand{\PET}{\operatorname{PET}}
\newcommand{\Supp}{\operatorname{Supp}}

\newcommand{\mult}{\operatorname{mult}}

\newcommand{\Bb}{\mathcal{B}}
\newcommand{\Ff}{\mathcal{F}}
\newcommand{\ff}{\mathfrak{f}}

\newcommand{\Ii}{\Gamma}

\newcommand{\Xx}{\mathcal{X}}
\newcommand{\Ss}{\mathcal{S}}

\newtheorem{thm}{Theorem}[section]

\newtheorem{cor}[thm]{Corollary}
\newtheorem{lem}[thm]{Lemma}
\newtheorem{prop}[thm]{Proposition}

\theoremstyle{definition}
\newtheorem{defn}[thm]{Definition}
\newtheorem{ques}[thm]{Question}
\theoremstyle{definition}
\newtheorem{rem}[thm]{Remark}

\newtheorem{defthm}[thm]{Definition-Theorem}

\newtheorem{nota}[thm]{Notation}

\theoremstyle{definition}

\begin{document}

\title{Volume of algebraically integrable foliations and locally stable families}
\author{Jingjun Han, Junpeng Jiao, Mengchu Li, Jihao Liu}

\address{Shanghai Center for Mathematical Sciences, Fudan University, Jiangwan Campus, Shanghai, 200438, China}
\email{hanjingjun@fudan.edu.cn}

\address{Yau Mathematical Sciences Center, Tsinghua University, Hai Dian District, Beijing 100084, China}
\email{jiao-jp@tsinghua.edu.cn}

\address{School of Mathematical Sciences, Fudan University, Shanghai 200433, China}
\email{mengchuli21@m.fudan.edu.cn}

\address{Department of Mathematics, Northwestern University, 2033 Sheridan Rd, Evanston, IL 60208, USA}
\email{jliu@northwestern.edu}

\subjclass[2020]{14E30, 37F75}
\keywords{Volume, Algebraically integrable foliation, Locally stable family, Moduli theory}
\date{\today}

\begin{abstract}
In this paper, we study the volume of algebraically integrable foliations and locally stable families. We show that, for any canonical algebraically integrable foliation, its volume belongs to a discrete set depending only on its rank and the volume of its general leaves. In particular, if the foliation is of general type, then its volume has a positive lower bound depending only on its rank and the volume of its general leaves. This implies some special cases of a question posed by Cascini, Hacon, and Langer.

As a consequence, we show that the relative volume of a stable family with a normal generic fiber belongs to a discrete set if the dimension and the volume of its general fibers are bounded. Log versions of the aforementioned theorems are also provided and proved.
\end{abstract}

\maketitle

\tableofcontents

\section{Introduction}
We work over the field of complex numbers $\mathbb{C}$.

The \emph{volume} of the canonical divisor of a projective variety $X$, or simply the volume of $X$, is an important invariant in birational geometry and the minimal model program. Hacon, M\textsuperscript{c}Kernan, and Xu famously proved that the volumes of log canonical varieties of fixed dimension satisfy the descending chain condition (DCC), and in particular, have a positive lower bound when the varieties are of general type \cite[Theorem 1.3]{HMX14}. These results are crucial for establishing the moduli theory for stable varieties.

In recent developments of the minimal model program and birational geometry, many new structures have been introduced, and it is natural to ask whether we can establish a good moduli theory for these structures. As a first step, we would like to investigate whether similar properties of volumes also hold for these new structures. For example, the volume of the anti-canonical divisor for Fano varieties has been proven to be well-controlled \cite{Bir19}, which is a crucial step in the proof of the BAB conjecture \cite{Bir19,Bir21a} and hence for the moduli theory of Fano varieties.

In this paper, we focus on the behavior of the volume of two related structures: algebraically integrable foliations and locally stable families.

\smallskip

\noindent\textbf{Volume of Algebraically Integrable Foliations}. A \emph{foliation} $\Ff$ on a normal projective variety $X$ is a coherent, saturated subsheaf of $T_X$ that is closed under the Lie bracket. Each foliation $\Ff$ is associated with a foliated canonical divisor $K_{\Ff}$, and the singularities of $\mathcal{F}$ can be defined similarly to the singularities of usual varieties \cite{McQ08} by considering $K_{\Ff}$ instead of $K_X$. Therefore, we can study the birational geometry of algebraically integrable foliations based on the behavior of $K_{\Ff}$ instead of $K_X$ (cf. \cite{ACSS21,CS20,CS21,CHLX23,LMX24,Spi20,SS22,SS23}).

Hacon and Langer asked whether the volume of any surface foliation with canonical singularities is well-ordered (e.g., satisfies the DCC) and has a lower bound \cite[Question 4]{HL21}. This question was generalized to higher dimensions by Cascini \cite[Section 3]{Cas21}.

\begin{ques}[{\cite[Section 3]{Cas21}; see also \cite[Question 4]{HL21}}]\label{ques: vol kf lower bound}
    Let $d$ be a positive integer and let $\Ff$ be a foliation with canonical singularities of general type (i.e. $K_{\Ff}$ is big) on a normal projective variety $X$ of dimension $d$.
    \begin{enumerate}
        \item Is there a positive real number $\epsilon=\epsilon(d)$ such that $\vol(K_{\Ff})\geq \epsilon$?
        \item Does $\vol(K_{\Ff})$ satisfy the DCC (also called well-ordered in \cite{HL21})?
    \end{enumerate}
\end{ques}

Question \ref{ques: vol kf lower bound} is very difficult and is not even known for rank one foliations on surfaces, although many works have been dedicated to it \cite{PS19,Che21,HL21,LT22,Lü24}.

In this paper, we consider a special case of Question \ref{ques: vol kf lower bound}, namely, the case when the volume of the general leaves of $\Ff$ is bounded. From the moduli theory point of view, this is a natural condition to ask, as the general leaves of any bounded family of foliations are always bounded. In this case, we not only provide a positive answer to Question \ref{ques: vol kf lower bound} but can further prove that $\vol(K_{\Ff})$ belongs to a discrete set. Of course, algebraic integrability is necessary in this case for us to define the ``volume of general leaves" as the leaves must be algebraic. We have the following result:

\begin{thm}\label{thm: main theorem foliation no boundary}
Let $\Ff$ be an algebraically integrable foliation with canonical singularities on a normal projective variety $X$. Let $L$ be a general leaf of $\Ff$. Then there exists a discrete set $\Ii_0$ depending only on $\rk\Ff$ and $\vol(K_L)$ such that $\vol(K_{\Ff})\in\Ii_0$.

In particular, $\vol(K_{\Ff})$ belongs to a DCC set, and there exists a positive real number $\delta$ depending only on $\rk\Ff$ and $\vol(K_L)$, such that $\vol(K_{\Ff})\geq\delta$ if $\Ff$ is of general type. 
\end{thm}

Theorem \ref{thm: main theorem foliation no boundary} answers Cascini-Hacon-Langer's Question \ref{ques: vol kf lower bound} for algebraically integrable foliations whose leaves have bounded volumes, and in particular, for algebraically integrable foliations with bounded leaves. Note that, in Theorem \ref{thm: main theorem foliation no boundary}, we only need to fix $\rk \Ff$ rather than $\dim X$.

\begin{rem}
    It seems that, if we do not have the algebraically integrable condition, Question \ref{ques: vol kf lower bound} will be more difficult to solve. This is due to several reasons. First, we do not know the minimal model program for non-algebraically integrable foliations in dimension $\geq 4$. Second, the general leaf $L$ may not be integrable, making it difficult to consider its canonical divisor $K_L$.
    
    We also expect more interesting behavior in the volume to appear for non-algebraically integrable foliations. For example, \cite[Page 207]{HL21} shows that there exists a sequence of foliations $\{\Ff_n\}_{n=2}^{\infty}$ with canonical singularities on surfaces, such that $K_{\Ff_n}$ is ample for each $n$, $\vol(K_{\Ff_n})$ is strictly increasing, and $\lim_{n\rightarrow+\infty}\vol(K_{\Ff_n})=1$. However, these foliations are quotients of Jouanolou’s foliations, which are known to be not algebraically integrable \cite{Jou79}.
\end{rem}

\smallskip
\noindent\textbf{Volume of Locally Stable Families}. The theory of moduli for stable varieties and pairs has been established over the past several decades. For a comprehensive discussion, we refer the reader to \cite{Kol23} and the references therein, and to \cite{Bir22} for a more general theory. A key structure in moduli theory is the notion of \emph{locally stable families}, which provides a nice parametrization of algebraic varieties belonging to a bounded family. For simplicity, in this paper, we only consider locally stable families over reduced bases (see Definition \ref{defn: lsf}).

Given a locally stable family $f: X \rightarrow Z$, the ambient variety $X$ may have bad singularities, i.e., singularities that are worse than log canonical (lc). However, many results in birational geometry still hold if we consider the relative canonical divisor $K_{X/Z}$. In particular, the existence of the minimal model program for locally stable families with a normal ambient variety was recently proven by Meng and Zhuang and has essentially been applied to wall-crossing theory \cite{MZ23}. Inspired by their work, we aim to study the birational geometry of the relative canonical divisor $K_{X/Z}$ in more detail.

One natural concept to consider is $\vol(K_{X/Z})$, the volume of $K_{X/Z}$. This is because $\vol(K_{X/Z})$ is preserved under any sequence of steps of $K_{X/Z}$-MMPs, as well as birational base changes. From a geometric standpoint, $\vol(K_{X/Z})$ measures how far the locally stable family is from being isotrivial. In particular, it is interesting to ask whether some numerical properties of the set of $\vol(K_X)$, such as the DCC property and the existence of a positive lower bound \cite[Theorem 1.3]{HMX14}, also hold for $\vol(K_{X/Z})$.

As a consequence of Theorem \ref{thm: main theorem foliation no boundary}, we show that $\vol(K_{X/Z})$ belongs to a discrete set under the extra condition that the general fiber of $f: X \rightarrow Z$ has fixed volume, and in particular, it satisfies the DCC. This is because the foliation induced by $X \rightarrow Z$ is an lc algebraically integrable foliation and $K_{\Ff}=K_{X/Z}$.

\begin{thm}\label{thm: main theorem lsf no boundary}
        Let $f: X\rightarrow Z$ be a stable family with normal generic fiber, such that $Z$ is a reduced projective scheme. Let $X_g$ be a general fiber of $f$. Then there exists a discrete set $\Ii_0$ depending only on $\dim X_g$ and $\vol(K_{X_g})$, such that $\vol(K_{X/Z})\in\Ii_0$. 

        In particular, $\vol(K_{X/Z})$ belongs to a DCC set, and there exists a positive real number $\delta$ depending only on $\dim X_g$ and $\vol(K_{X_g})$, such that $\vol(K_{X/Z})\geq\delta$ if $f: X\rightarrow Z$ is of maximal variation.
\end{thm}
Theorem \ref{thm: main theorem lsf no boundary} can be viewed as an analogue of \cite[Theorem 1.3]{HMX14} for stable families. We remark that the slope formulas obtained in \cite[Theorems A,B]{CTV23} share some similarities with the inequality we obtained in Theorem \ref{thm: main theorem lsf no boundary}. Compared to our result, \cite[Theorems A,B]{CTV23} do not require the boundedness of the general fiber but require $Z$ to be a curve, and the inequalities rely on $\deg(f_*\mathcal{O}_X(K_{X/Z}))$ (or $\frac{1}{m^{\dim X}} \cdot \deg(f_*\mathcal{O}_X(mK_{X/Z}))$).

\medskip

\noindent\textbf{Log Versions}. Instead of considering varieties and foliations only, in practice, it is sometimes more natural to consider their ``log versions," that is, pairs and foliated triples, respectively. We would like to ask whether the log version of Theorems \ref{thm: main theorem lsf no boundary} and \ref{thm: main theorem foliation no boundary} also hold. The answer is positive: we have the following foliated triple version of Theorem \ref{thm: main theorem foliation no boundary}. Moreover, we can relax the condition ``canonical" in Theorem \ref{thm: main theorem foliation no boundary} to ``lc" (see Definition \ref{defn: foliation singularity}).

\begin{thm}\label{thm: main theorem foliation with boundary}
Let $d$ be a positive integer, $\Ii\subset [0,1]$ a DCC set, and $v$ a positive real number. Then there exists a positive real number $\delta$ and a discrete set $\Ii_0$ depending only on $d,\Ii$ and $v$ satisfying the following. Assume that 
\begin{itemize}
    \item $(X,\Ff,B)$ is a projective lc algebraically integrable foliated triple,
    \item  $L$ is the normalization of a general leaf of $\Ff$ and $K_L+B_L:=(K_{\Ff}+B)|_L$,
    \item $\rk\Ff=d$ and the coefficients of $B$ belong to $\Ii$, and
    \item $\vol(K_L+B_L)=v$.
\end{itemize}
Then:
\begin{enumerate}
    \item If $K_{\Ff}+B$ is big, then $\vol(K_{\Ff}+B)\geq\delta$.
    \item If $(X_L,B_L)$ has a good minimal model, then $\vol(K_{\Ff}+B)\in\Ii_0$.
\end{enumerate}
\end{thm}
We remark that, since $K_L + B_L$ is big, the condition ``$(X_L, B_L)$ has a good minimal model” in condition (2) always holds either when $(X_L, B_L)$ is klt, or when $K_L + B_L$ is already ample.

As a direct consequence of Theorem \ref{thm: main theorem foliation with boundary}, we obtain the following log version of Theorem \ref{thm: main theorem lsf no boundary}:

\begin{thm}\label{thm: main theorem lsf with boundary}
Let $d$ be a positive integer, $\Ii\subset [0,1]$ a DCC set, and $v$ a positive real number. Then there exists a positive real number $\delta$ and a discrete set $\Ii_0$ depending only on $d,\Ii$ and $v$ satisfying the following. Assume that 
\begin{itemize}
\item $f: (X,B)\rightarrow Z$ be a locally stable family with normal generic fiber,
\item $Z$ is a reduced projective scheme and $(X_g,B_g)$ is a general fiber of $f: (X,B)\rightarrow Z$,
\item $\dim X_g=d$ and the coefficients of $B$ belong to $\Ii$, and
\item $\vol(K_{X_g}+B_g)=v$.
\end{itemize}
Then:
\begin{enumerate}
    \item If $K_{X/Z}+B$ is big, then $\vol(K_{X/Z}+B)\geq\delta$.
    \item If $(X_g,B_g)$ has a good minimal model, then $\vol(K_{X/Z}+B)\in\Ii_0$.
\end{enumerate}
\end{thm}
We remark that, since $K_{X_g} + B_g$ is big, the condition ``$(X_g, B_g)$ has a good minimal model” in condition (2) always holds either when $(X_g, B_g)$ is klt or when $K_{X_g} + B_g$ is already ample (e.g., when $f: (X, B) \rightarrow Z$ is stable).

\medskip

\noindent\textit{Idea of the Proofs}. We provide a sketch of the proof of Theorem \ref{thm: main theorem foliation with boundary}. Theorems \ref{thm: main theorem lsf no boundary}, \ref{thm: main theorem foliation no boundary}, and \ref{thm: main theorem lsf with boundary} are direct consequences of Theorem \ref{thm: main theorem foliation with boundary} and the minimal model program for algebraically integrable foliations.

First, we may reduce to the case when $X$ is $\Qq$-factorial klt by taking an ACSS modification \cite{ACSS21,CHLX23}, and $\Ff$ is induced by a contraction $f: X \rightarrow Z$. By the ACC for lc thresholds for algebraically integrable foliations \cite{DLM23}, we can prove the ACC for pseudo-effective thresholds for algebraically integrable foliations (Theorem \ref{thm: acc pet aif}), and reduce (1) to the case when $(X_L, B_L)$ is klt. Now (1) follows from (2), so we only need to prove (2).

Thanks to the establishment of the minimal model program for algebraically integrable foliations \cite{ACSS21,CS23,CHLX23,LMX24}, we are able to show that, when $(X_L, B_L)$ has a good minimal model, $(X, \Ff, B)$ has a minimal model $(X', \Ff', B')$ that is a good minimal model over $Z$. Since volume is preserved under the MMP, we may replace $(X, \Ff, B)$ with the ample model of $(X', \Ff', B')$ over $Z$ and assume that $K_{\Ff} + B$ is ample over $Z$. In this case, generically over $Z$, $f: (X, B) \rightarrow Z$ is a stable family. More precisely, we reduce Theorem \ref{thm: main theorem foliation with boundary} to the special case, Theorem \ref{thm: main theorem foliation already induced case}.

Our next goal is to use moduli theory for stable varieties and the lc property of the foliation to reduce to the case when $f: (X, B) \rightarrow Z$ is a stable family and $Z \subset \Ss$, where $\Ss \rightarrow M$ is a finite cover of bounded degree, and $M$ is the coarse moduli space of the functor that parametrizes the general log leaves of $(X, \Ff, B)$. The idea is in two steps.

First, we take a resolution of $Z' \rightarrow Z$ so that the induced map $Z' \rightarrow M$ is a morphism. Then we consider $\bar{Z} := Z' \times_{M} \Ss$, and consider the foliation induced by $\bar{X} \rightarrow \bar{Z}$, where $\bar{X}$ is the normalization of the main component of $X \times_Z Z'$. A key observation at this point is that both the base change $Z' \rightarrow Z$ and the base change $\bar{Z} \rightarrow Z'$ will preserve the log canonical property of $(X, \Ff, B)$, and moreover, the boundary of the foliated triple still has coefficients in $\Ii$. More precisely, for the base change $Z' \rightarrow Z$, which is a resolution, we may apply \cite[Lemma 3.5]{LMX24} (see Lemma \ref{lem: ACSS exist core mod equidim}), and for the base change $\bar{Z} \rightarrow Z'$, which is a finite cover, we may apply \cite[Lemma 3.4, 4.3]{Dru21} (see Proposition \ref{prop: cover formula}). Moreover, the degree of $\bar{Z} \rightarrow Z$ is equal to the degree of $\Ss \rightarrow M$, which is a constant. Therefore, we are free to replace $X$ with $\bar{X}$ and assume that there exists a morphism $Z \rightarrow \Ss$. After a further birational base change, we may assume that $Z$ is smooth. Moreover, there exists an induced stable family $(\Xx_Z, \Bb_Z) \rightarrow Z$ which is the base change of the universal family $(\Xx, \Bb) \rightarrow \Ss$.

Second, we can take a finite cover of $Z$ along the image of the ramification divisor $R(f)$ on $Z$ and obtain a stable family $(X', B') \rightarrow Z'$. This again follows from \cite[Lemma 3.4, 4.3]{Dru21} and the fact that $(X, \Ff, B)$ is lc, and was rigorously stated in \cite[Proposition 3.5]{FS22} (see Proposition \ref{prop: special cover}). By \cite[Proposition 8.64]{Kol23}, after another cover $\tilde{Z} \rightarrow Z'$, we may obtain a model $(\tilde{X}, \tilde{B}) \rightarrow \tilde{Z}$ that is isomorphic to the corresponding base change of $(\Xx_Z, \Bb_Z) \rightarrow Z$. An easy computation indicates that $\vol(K_{\Xx_Z/Z} + \Bb_Z) = \vol(K_{\Ff} + B)$ (see Lemma \ref{lem: over same base lsf volume are same}) so we may replace $(X, B) \rightarrow Z$ with $(\Xx_Z, \Bb_Z) \rightarrow Z$ and assume that $(X, B) \rightarrow Z$ is a stable family and is the base change of $(\Xx, \Bb) \rightarrow \Ss$.

Finally, we let $T$ be the image of $Z$ in $\Ss$. Then $Z \rightarrow T$ is a finite cover, so $l \cdot \vol(K_{X_T/T} + B_{T}) = \vol(K_X + B)$ for some integer $l > 0$, where $(X_{T}, B_{T}) \rightarrow T$ is the induced stable family. Now we only need to prove the discreteness of $\vol(K_{X_T/T} + B_{T})$. If $\Bb$ is a $\Qq$-divisor, then the Cartier index of $K_{X_T/T} + B_T$ is bounded by the Cartier index of $K_{\Xx/\Ss} + \Bb$ and we are done. When $\Bb$ has irrational coefficients, we may perturb the coefficients of $\Bb$ and get stable families over $\Ss$ with rational coefficients (which follows from the openness of ampleness and the existence of uniform lc rational polytopes \cite[Theorem 5.6]{HLS19}). The discreteness of $\vol(K_{X_T/T} + B_{T})$ now follows from an easy intersection number computation.

\begin{rem}
After the first draft of this paper, L. Tasin informed us that G. Codogni, Z. Patakfalvi, and himself have obtained the existence of a positive lower bound for $\vol(K_{X/Z} + B)$ for stable families $f: (X, B) \rightarrow Z$, provided that the general fiber $(X_g, B_g)$ of $f: (X, B) \rightarrow Z$ is klt and $\Ii \subset \Qq$ \cite{CPT24}. The methods used in their proof differ from our proof of Theorem \ref{thm: main theorem lsf with boundary}.
\end{rem}

\medskip

\noindent\textbf{Acknowledgement}. The authors would like to thank Paolo Cascini, Guodu Chen, Chen Jiang, Yuchen Liu, Fanjun Meng, Calum Spicer, Roberto Svaldi, Luca Tasin, Lingyao Xie, Qingyuan Xue, and Ziquan Zhuang for useful discussions. The third author completed this work as a visiting student at Imperial College London. He expresses his gratitude to his advisor Professor Meng Chen for his great support and encouragement, and to Professor Paolo Cascini for the numerous helpful discussions and support. The first author is supported by NSFC for Excellent Young Scientists (\#12322102), and the National Key Research and Development Program of China (\#2023YFA1010600, \#2020YFA0713200). The second author was supported by grants from Tsinghua University, Yau Mathematical Science Center. The third author is supported by the China Scholarship Council (No. 202306100155) and the Fudan Elite PhD Program.

\section{Locally stable families}

We will adopt the standard notations and definitions on MMP in \cite{KM98,BCHM10} and use them freely. We refer the reader to \cite{Kol23} for basic notations and properties of moduli theory but we shall recall some important concepts here.

\subsection{Preliminaries}

\begin{defn}
    A \emph{contraction} is a projective morphism $f: X\rightarrow Y$ such that $f_*\mathcal{O}_X=\mathcal{O}_Y$. A \emph{finite cover} is a finite surjective morphism. 
\end{defn}

\begin{nota}
    Let $f: X\dashrightarrow X'$ be a birational map between normal schemes. We denote by $\Exc(f)$ the reduced divisor supported on the codimension one part of the exceptional locus of $f$.
\end{nota}

\begin{nota}
    Let $X$ be a reduced scheme and $D,D'$ two $\Rr$-divisors on $X$. We define 
    $D\wedge D':=\sum_P\min\{\mult_PD,\mult_PD'\}P$ where the sum runs through all the prime divisors $P$ on $X$. We denote by $\Supp D$ the reduced divisor supported on $D$.
\end{nota}

\begin{defn}
Let $m$ be a positive integer and $\bm{v}\in\mathbb R^m$. The \emph{rational envelope} of $\bm{v}$ is the minimal rational affine subspace of $\mathbb R^m$ which contains $\bm{v}$. For example, if $m=2$ and $\bm{v}=\left(\frac{\sqrt{2}}{2},1-\frac{\sqrt{2}}{2}\right)$, then the rational envelope of $\bm{v}$ is $(x_1+x_2=1)\subset\mathbb R^2_{x_1x_2}$.
\end{defn}

\begin{defn}
Let $X$ be a projective normal variety of dimension $d$, and $D$ an $\Rr$-divisor on $X$. The \emph{volume} of $D$ is defined as: $$\vol(D):=\limsup_{m\to \infty} \frac{h^0(\lfloor mD\rfloor)}{m^d/d!}
.$$
By \cite[Theorem 3.5]{FKL16}, we can take the limit superior in the above definition as the limit, and when $D$ is $\Rr$-Cartier, the above definition agrees with the definition in \cite[Corollary 2.4.5]{Laz04}. 
\end{defn}

\begin{lem}
If $f: X \to Y$ is a generically finite proper morphism from a normal projective variety $X$, and $D$ an $\Rr$-Cartier $\Rr$-divsior on $Y$, then 
$$\vol(f^{*}D)=\deg f\cdot \vol (D).$$
\end{lem}

\begin{proof}
This follows from the homogeneity of asymptotic cohomological functions \cite[Proposition 2.7]{Kur06}, and asymptotic cohomological functions of pullbacks \cite[Proposition 2.9]{Kur06}.
\end{proof}

\subsection{Locally stable families}

\begin{defn}[Slc pair]\label{defn: slc pair}
In this paper, a \emph{pair} $(X,B)$ consists of a demi-normal scheme $X$ and an $\Rr$-divisor $B\geq 0$ such that $K_X+B$ is $\Rr$-Cartier. Let $\nu: X^\nu\rightarrow X$ be the normalization of $X$ and let $D\subset X$, $D^\nu\subset X^\nu$ be the conductors. We say that $(X,B)$ is \emph{slc} if 
\begin{enumerate}
\item $K_X+B$ is $\Rr$-Cartier,
    \item $\Supp B$ does not contain any irreducible component of the conductor $D$.
    \item $(X^\nu,B^\nu+D^\nu)$ is lc, where $B^\nu$ is the divisorial part of $\nu^{-1}(B)$.
\end{enumerate}
We say that $(X,B)$ is \emph{stable} if $(X,B)$ is projective, slc, and $K_X+B$ is ample. 
\end{defn}

\begin{defn}[Mumford divisor]\label{defn: mumford divisor}
Let $d$ be a positive integer and let $f: X\rightarrow Z$ be a flat morphism over a reduced scheme such that the fibers of $f$ are reduced, connected, $S_2$, and of pure dimension $d$. A \emph{Mumford divisor}$/Z$ is a divisor $D$ on $X$ satisfying the following.
\begin{enumerate}
    \item (Equidimensionality) Every irreducible component of $\Supp D$ dominates an irreducible component of $Z$, and all nonempty fibers of the induced morphism $f|_D: \Supp D\rightarrow Z$ are of pure dimension $d-1$.
    \item (Mumford) $f$ is smooth near the generic point of $f^{-1}(z)\cap\Supp D$ for any point $z\in Z$.
    \item (Generic Cartier) $D$ is Cartier near the generic points of $f^{-1}(z)\cap\Supp D$ for any point $z\in Z$.
\end{enumerate}
By \cite[Theorem 4.4]{Kol23}, if $Z$ is normal, then (1) and (2) imply (3).

An $\Rr$-divisor $B$ on $X$ is called a \emph{Mumford $\mathbb R$-divisor$/Z$} if $B=\sum a_iD_i$, where each $a_i\geq 0$ and each $D_i$ is a Mumford divisor$/Z$.
\end{defn}

\begin{defn}[Locally stable family]\label{defn: lsf}
A \emph{locally stable family} $f: (X,B)\rightarrow Z$ consists of a flat projective morphism $ f: X\rightarrow Z$ over a reduced scheme with demi-normal and connected fibers, and a Mumford $\Rr$-divisor$/Z$ $B$ on $X$, such that
\begin{enumerate}
\item $K_{X/Z}+B$ is $\Rr$-Cartier, and
\item $(X_z,B_z)$ is slc for any $z\in Z$, where $X_z=f^{-1}(z)$ and $B_z=B|_{X_z}$ (we refer the reader to \cite[Section 4.1]{Kol23} for the definition of such restriction).
\end{enumerate}
Here $K_{X/Z}$ is the relative canonical divisor corresponding to the relative dualizing sheaf $\omega_{X/Z}$ and we refer the reader to \cite[2.68]{Kol23} for more details.

A \emph{stable family} is a locally stable family $f: (X,B)\rightarrow Z$ such that $K_{X/Z}+B$ is ample$/Z$.
\end{defn}

\subsection{Moduli of stable families}

\begin{defn}
    Let $\Ii\subset [0,1]$ be a set. We define
    $$\Ii_+:=\{0\}\cup\left(\left\{\sum_{i=1}^n\gamma_i \mid n\in\mathbb N^+,\gamma_1,\dots,\gamma_n\in\Ii\right\}\cap [0,1]\right).$$
\end{defn}

\begin{defn}[$(d,\Ii,v)$-stable pairs]
    Let $d$ be a positive integer, $\Ii\subset [0,1]$ a set such that $\Ii=\Ii_+$, $v$ a positive real number, and $(X,B)$ a pair. We say that $(X,B)$ is a \emph{$(d,\Ii,v)$-stable pair} if  $\dim X=d$, $(X,B)$ is stable, the coefficients of $B$ belong to $\Ii$, and  $\vol(K_X+B)=v$. A \emph{family of $(d,\Ii,v)$-stable pairs} is a locally stable family $f: (X,B)\rightarrow Z$, such that
 \begin{enumerate}
     \item $B=\sum b_iB_i$ where each $b_i\in\Ii$ and each $B_i$ is a Mumford divisor$/Z$, and
     \item $(X_z,B_z)$ is a $(d,\Ii,v)$-stable pair for any point $z\in Z$, where $X_z=f^{-1}(z)$ and $B_z=B|_{X_z}$.
 \end{enumerate}
We denote by $\mathcal{S}(d,\Ii,v)$ the set of all $(d,\Ii,v)$-stable pairs. We denote by $\mathfrak{S}(d,\Ii,v)$ the moduli functor by setting
 $$\mathfrak{S}(d,\Ii,v)(Z)=\{\text{families of }(d,\Ii,v)\text{-stable pairs over }Z\}$$
 where $Z$ is any reduced scheme. 
\end{defn}

\begin{defthm}\label{defthm: kol23 4.1}
Let $d$ be a positive integer, $\Ii\subset [0,1]$ a DCC set such that $\Ii=\Ii_+$, and $v$ a positive real number. Then the functor $\mathfrak{S}(d,\Ii,v)$ has a projective coarse moduli space $M(d,\Ii,v)$. For any family $f: (X,B)\rightarrow Z$ of $(d,\Ii,v)$-stable pairs, we call the induced morphism $Z\rightarrow M(d,\Ii,v)$ the corresponding \emph{moduli map}.
\end{defthm}
\begin{proof}
By \cite[6.8.4]{Kol23} (the real coefficient case of \cite[Theorem 1.1]{HMX18}; see also \cite[Theorem 6.4]{Li24}), there exists a finite set $\Ii_0\subset\Ii$ depending only on $d,\Ii$ and $v$ such that any $(d,\Ii,v)$-stable family is a $(d,\Ii_0,v)$-stable family. Now the existence of the projective coarse moduli space follows from \cite[Theorem 4.1]{Kol23}.
\end{proof}

\begin{thm}[{\cite[Corollary 6.19]{KP17}}]\label{kp17 6.19}
Let $d$ be a positive integer, $\Ii\subset [0,1]$ a DCC set such that $\Ii=\Ii_+$, and $v$ a positive real number. Then there exists a reduced scheme $Z$ and $f\in\mathfrak{S}(d,\Ii,v)(Z)$, such that the corresponding moduli map $Z\rightarrow M(d,\Ii,v)$ is a finite cover.
\end{thm}
\begin{proof}
It is essentially \cite[Corollary 6.19]{KP17}. We remark that \cite[Corollary 6.19]{KP17} requires the stable objects we parametrize to have rational coefficients, and require that the functor to parametrizing normal bases rather than reduced scheme bases. However the real coefficient and the reduced scheme base case follows from the same lines of the proof by using the fact that there exists a good moduli theory for $(d,\Ii,v)$-stable pairs (Definition-Theorem \ref{defthm: kol23 4.1}).
\end{proof}

\subsection{Base change}

\begin{defn}[Base change]
    Let $f: (X,B)\rightarrow Z$ be a (locally) stable family. A \emph{base change} of $f: (X,B)\rightarrow Z$ is $f': (X',B')\rightarrow Z'$ satisfying the following:
    \begin{itemize}
        \item There is a morphism $h_Z: Z'\rightarrow Z$ such that $X'=X\times_ZZ'$, and $h: X'\rightarrow X$ is the induced morphism.
        \item $B'$ is the unique $\Rr$-divisor such that $K_{X'/Z'}+B'=h^*(K_{X/Z}+B)$.
    \end{itemize}
By definition (cf. \cite[4.6.2]{Kol23}), $f': (X',B')\rightarrow Z'$ is also a (locally) stable family.  If $\dim Z'=\dim Z$ and $h_Z$ is surjective, then we say that $f': (X',B')\rightarrow Z'$ is a \emph{birational base change} of $f: (X,B)\rightarrow Z$.
\end{defn}

\begin{prop}[{\cite[Proposition 8.64]{Kol23}}]\label{prop: isom functor is finite}
Let $d$ be a positive integer, $\Ii\subset [0,1]$ a DCC set such that $\Ii=\Ii_+$, and $v$ a positive real number. Let $f_i: (X_i,B_i)\rightarrow Z$, $i\in\{1,2\}$ be two families of  $(d,\Ii,v)$-stable pairs $f_i:(X_i,B_i)\rightarrow Z$ with corresponding moduli maps $g_i: Z\rightarrow M(d,\Ii,v)$, such that $g_1=g_2$. Then there exists a finite cover $\pi: \tilde Z\rightarrow Z$, such that the induced base changes
$$\tilde f_i: (X_i,B_i)\times_{Z}\tilde Z\rightarrow\tilde Z$$
are isomorphic. In particular, the generic fiber of $f_1$ is normal if and only if the generic fiber of $f_2$ is normal.
\end{prop}

\begin{lem}\label{lem: over same base lsf volume are same}
Let $d$ be a positive integer, $\Ii\subset [0,1]$ a DCC set such that $\Ii=\Ii_+$, and $v$ a positive real number. Let $f_i: (X_i,B_i)\rightarrow Z$, $i\in\{1,2\}$ be two families of  $(d,\Ii,v)$-stable pairs $f_i:(X_i,B_i)\rightarrow Z$ with corresponding moduli maps $g_i: Z\rightarrow M(d,\Ii,v)$, such that $g_1=g_2$. Then
$$\vol(K_{X_1/Z}+B_1)=\vol(K_{X_2/Z}+B_2).$$
\end{lem}
\begin{proof}
By Proposition \ref{prop: isom functor is finite}, there exists a finite cover $\pi: \tilde Z\rightarrow Z$, such that the induced base changes
$$\tilde f_i: (\tilde X_i,\tilde B_i):=(X_i,B_i)\times_{Z}\tilde Z\rightarrow\tilde Z$$
are isomorphic.
Then
$$\vol(K_{X_1/Z}+B_1)=\frac{1}{\deg\pi}\vol(K_{\tilde X_1/\tilde Z}+\tilde B_1)=\frac{1}{\deg\pi}\vol(K_{\tilde X_2/\tilde Z}+\tilde B_2)=\vol(K_{X_2/Z}+B_2).$$
\end{proof}

\subsection{Maximal variation and bigness}

\begin{defn}
Let $f: (X,B)\rightarrow Z$ be a stable family. The \emph{variation} of $f: (X,B)\rightarrow Z$ is $\dim Z-d$ where $d$ is the dimension of a general isomorphism equivalence class of the log fibers. We say that $f: (X,B)\rightarrow Z$ is \emph{of maximal variation} if the variation of $f: (X,B)\rightarrow Z$ equals to $\dim Z$.
\end{defn}

\begin{prop}[cf. {\cite[Proposition 2.15]{PX17}}]\label{prop: maximal variation big and nef}
Let $f: (X,B)\rightarrow Z$ be a stable family of maximal variation such that 
\begin{enumerate}
    \item $Z$ is projective, and
    \item the generic fiber $(X_{\eta},B_{\eta})$ of $f: (X,B)\rightarrow Z$ is lc.
\end{enumerate}
Then $K_{X/Z}+B$ is big and nef.
\end{prop}
\begin{proof}Since the bigness and nefness of $K_{X/Z}+B$ is preserved under birational base change of $f: (X,B)\rightarrow Z$, possibly replacing $Z$ with a normalization followed by a resolution, and replace $f: (X,B)\rightarrow Z$ with the corresponding base change, we may assume that $Z$ is smooth projective.

Let $B=\sum_{i=1}^m b_iB_i$ where $B_i$ are the irreducible components of $B$ and let $B(\bm{v})=\sum v_iB_i$ for any $\bm{v}=(v_1,\dots,v_m)\in\mathbb R^m$. By \cite[Lemma 5.3]{HLS19}, $K_{X/Z}+B(\bm{v})$ is $\Rr$-Cartier for any $\bm{v}$ in the rational envelope of $\bm{b}:=(b_1,\dots,b_m)$. Since ampleness is an open condition, there exists an open subset $U\ni\bm{b}$ of the rational envelope of $\bm{b}$ such that $K_{X/Z}+B(\bm{v})$ is ample$/Z$ for any $\bm{v}\in U$. Since $Z$ is smooth, for any snc divisor $D$ on $Z$, $(X,B+f^*D)$ is slc. By \cite[Theorem 5.6]{HLS19}, possibly shrinking $U$, we may assume that $(X,B(\bm{v})+f^*D)$ is slc for any $\bm{v}\in U$ and any snc divisor $D$ on $Z$. Therefore, $f: (X,B(\bm{v}))\rightarrow Z$ is a locally stable family. Pick vectors $\bm{v}_1,\dots,\bm{v}_{m+1}\in U\cap\mathbb Q^m$ such that $\bm{b}$ is in the interior of the convex hull of $\bm{v}_1,\dots,\bm{v}_{m+1}$. By \cite[Proposition 2.15]{PX17}, $K_{X/Z}+B(\bm{v}_i)$ is big and nef for any $i$. Thus $K_{X/Z}+B$ is big and nef.
\end{proof}

\begin{lem}\label{lem: non maximal variation not big}
    Let $f: (X,B)\rightarrow Z$ be a stable family that is not of maximal variation such that $Z$ is normal projective. Then $K_{X/Z}+B$ is not big.
\end{lem}
\begin{proof}
Since bigness of $K_{X/Z}+B$ is preserved under a finite base change, by \cite[Corollary 6.20]{KP17}, possibly passing to a finite base change, we may assume that there are two morphisms $h: X\rightarrow X'$ and $h_Z: Z\rightarrow Z'$ and a stable family $f': (X',B')\rightarrow Z'$ of maximal variation, such that $\dim Z'=\Var(f)<\dim Z$ and $f: (X,B)\rightarrow Z$ is a base change of  $f': (X',B')\rightarrow Z'$. Thus $\dim X>\dim X'$ and $h^*(K_{X'/Z'}+B')=K_{X/Z}+B$, so $K_{X/Z}+B$ is not big.
\end{proof}

\section{Algebraically integrable foliations}

\subsection{Basic definitions}

\begin{defn}[Foliations, {cf. \cite{ACSS21,CS21}}]\label{defn: foliation}
Let $X$ be a normal variety. A \emph{foliation} on $X$ is a coherent sheaf $\Ff\subset T_X$ such that
\begin{enumerate}
    \item $\Ff$ is saturated in $T_X$, i.e. $T_X/\Ff$ is torsion free, and
    \item $\Ff$ is closed under the Lie bracket.
\end{enumerate}
The \emph{rank} of the foliation $\Ff$ is the rank of $\Ff$ as a sheaf and is denoted by $\rk\Ff$. The \emph{canonical divisor} of $\Ff$ is a divisor $K_\Ff$ such that $\mathcal{O}_X(-K_{\mathcal{F}})\cong\mathrm{det}(\Ff)$. If $\Ff=0$, then we say that $\Ff$ is a \emph{foliation by points}.

Given any dominant map $h: Y\dashrightarrow X$, we denote by $h^{-1}\Ff$ the \emph{pullback} of $\Ff$ on $Y$ as constructed in \cite[3.2]{Dru21} and say that $h^{-1}\Ff$ is \emph{induced by} $\Ff$. Given any birational map $g: X\dashrightarrow X'$, we denote by $g_*\Ff:=(g^{-1})^{-1}\Ff$ the \emph{pushforward} of $\Ff$ on $X'$ and also say that $g_*\Ff$ is \emph{induced by} $\Ff$. We say that $\Ff$ is an \emph{algebraically integrable foliation} if there exists a dominant map $f: X\dashrightarrow Z$ such that $\Ff=f^{-1}\Ff_Z$, where $\Ff_Z$ is the foliation by points on $Z$, and we say that $\Ff$ is \emph{induced by} $f$.

A subvariety $S\subset X$ is called \emph{$\Ff$-invariant} if for any open subset $U\subset X$ and any section $\partial\in H^0(U,\Ff)$, we have $\partial(\mathcal{I}_{S\cap U})\subset \mathcal{I}_{S\cap U}$, 
where $\mathcal{I}_{S\cap U}$ is the ideal sheaf of $S\cap U$. For any prime divisor $P$ on $X$, we denote $\epsilon_{\Ff}(P):=1$ if $P$ is not $\Ff$-invariant and $\epsilon_{\Ff}(P):=0$ if $P$ is $\Ff$-invariant. For any prime divisor $E$ over $X$, we define $\epsilon_{\Ff}(E):=\epsilon_{\Ff_Y}(E)$ where $h: Y\dashrightarrow X$ is a birational map such that $E$ is on $Y$ and $\Ff_Y:=h^{-1}\Ff$.
\end{defn}

\begin{rem}
 Let $f: (X,B)\rightarrow Z$ be a (locally) stable family, $\Ff$ is induced by $f$, and $Z$ is normal. Then $K_{X/Z}\sim K_{\Ff}$ by \cite[Notation 2.7 and \S 2.9]{Dru17} and Definition \ref{defn: lsf}. In the rest of this paper, we will use this fact several times.
\end{rem}

\begin{defn}[{Tangent, {cf. \cite[Section 3.4]{ACSS21}}}]\label{defn: tangent to foliation}
 Let $X$ be a normal variety, $\Ff$ a foliation on $X$, and $V\subset X$ a subvariety. Suppose that $\Ff$ is a foliation induced by a dominant rational map $X\dashrightarrow Z$. We say that $V$ is \emph{tangent} to $\Ff$ if there exists a birational morphism $\mu: X'\rightarrow X$, an equidimensional contraction $f': X'\rightarrow Z'$, and a subvariety $V'\subset X'$, such that
    \begin{enumerate}
    \item $\mu^{-1}\Ff$ is induced by $f'$, and
        \item $V'$ is contained in a fiber of $f'$ and $\mu(V')=V$.
    \end{enumerate}
\end{defn}

\begin{defn}[Foliated triples]\label{defn: foliated triple}
 A \emph{foliated triple} $(X,\Ff,B)/U$ consists of a normal quasi-projective variety $X$, a foliation $\Ff$ on $X$, an $\Rr$-divisor $B\geq 0$ on $X$, and a projective morphism $X\rightarrow U$, such that $K_{\Ff}+B$ is $\mathbb R$-Cartier. If $U$ is not important, then we may drop $U$. 
 
 If $\Ff$ is algebraically integrable, then we say that $(X,\Ff,B)$ is algebraically integrable. If $X$ is $\Qq$-factorial, then we say that $(X,\Ff,B)$ is $\Qq$-factorial. If we allow $B$ to have negative coefficients, then we shall add the prefix ``sub-". If $B$ is a $\Qq$-divisor then we may add the prefix ``$\mathbb Q$-".
\end{defn}

\begin{defn}[Singularities]\label{defn: foliation singularity}
Let $(X,\Ff,B)$ be a foliated triple. For any prime divisor $E$ over $X$, let $f: Y\rightarrow X$ be a birational morphism such that $E$ is on $Y$, and suppose that
$$K_{\Ff_Y}+B_Y:=f^*(K_\Ff+B)$$
where $\Ff_Y:=f^{-1}\Ff$. We define $a(E,\Ff,B):=-\mult_EB_Y$ to be the \emph{discrepancy} of $E$ with respect to $(X,\Ff,B)$. If $\Ff=T_X$, then we define $a(E,X,B):=a(E,\Ff,B)$ which is the usual discrepancy for pairs.

We say that $(X,\Ff,B)$ is \emph{lc} (resp. \emph{klt}) if $a(E,\Ff,B)\geq -1$ (resp. $>-1$) for any prime divisor $E$ over $X$.  $(X,\Ff,B)$ is \emph{canonical} (resp. \emph{terminal}) if $a(E,\Ff,B)\geq 0$ (resp. $>0$) for any prime divisor $E$ that is exceptional over $X$.  For foliated sub-triples, we define singularities in the same way and we shall add the prefix ``sub-" for the descriptions of singularities. Note that our definition of klt aligns with \cite{LMX24} and is different from other references. 

An \emph{lc place} of $(X,\Ff,B)$ is a prime divisor $E$ over $X$ such that $a(E,\Ff,B)=-\epsilon_{\Ff}(E)$. An \emph{lc center} of $(X,\Ff,B)$ is the center of an lc place of $(X,\Ff,B)$ on $X$.
\end{defn}

\begin{lem}\label{lem: equivalent definition lc and klt}
Let $(X,\Ff,B)$ be a foliated sub-triple. Then:
\begin{enumerate}
 \item  $(X,\Ff,B)$ is sub-lc if and only if $a(E,\Ff,B)\geq-\epsilon_{\Ff}(E)$ for any prime divisor $E$ over $X$.
\item $(X,\Ff,B)$ is sub-klt if and only if $a(E,\Ff,B)>-1$ for any non-invariant prime divisor $E$ over $X$, and  $a(E,\Ff,B)\geq 0$ for any invariant prime divisor $E$ over $X$.
\end{enumerate}
\end{lem}
\begin{proof}
(1) is \cite[Lemma 2.10]{LMX24}. The if part of (2) folows from the definition so we only need to prove its only if part. Suppose the only if part of (2) does not hold, then there exists a sub-klt foliated sub-triple $(X,\Ff,B)$ such that $a(E,\Ff,B)<0$ for some $\Ff$-invariant prime divisor $E$ over $X$. By (1), $(X,\Ff,B)$ is not sub-lc, so it is not sub-klt, a contradiction.
\end{proof}

\subsection{Simple, core, and ACC modifications}

\begin{defn}[cf. Property $(*)$ foliations, {\cite[Definition 3.5]{ACSS21}, \cite[Definition 7.2.2]{CHLX23}}]\label{defn: foliation property *}
Let $(X,\Ff,B)$ be a foliated triple. Let $G\geq 0$ be a reduced divisor on $X$ and $f: X\rightarrow Z$ a contraction. We say that $(X,\Ff,B;G)/Z$ satisfies \emph{Property $(*)$} if the following conditions hold.
\begin{enumerate}
\item $\Ff$ is induced by $f$ and $G$ is an $\Ff$-invariant divisor.
\item $f(G)$ is of pure codimension $1$, $(Z,f(G))$ is log smooth, and $G=f^{-1}(f(G))$.
\item For any closed point $z\in Z$ and any reduced divisor  $\Sigma\ge f(G)$ on $Z$ such that  $(Z,\Sigma)$ is log smooth near $z$, $(X,B+G+f^*(\Sigma-f(G)))$ is lc over a neighborhood of $z$.
\end{enumerate}
We say that $f$, $Z$, and $G$ are \emph{associated} with $(X,\Ff,B)$.
\end{defn}

\begin{prop}[cf. {\cite[Proposition 3.6]{ACSS21}, \cite[Proposition 7.3.6]{CHLX23}}]\label{prop: weak cbf gfq}
Let $(X,\Ff,B)$ be a foliated triple. Let  $G\geq 0$ be a reduced divisor on $X$ and $f: X\rightarrow Z$ an equidimensional contraction, such that $(X,\Ff,B;G)/Z$ satisfies Property $(*)$ and $B$ is horizontal$/Z$. Then
$$K_{\Ff}+B\sim_{Z}K_X+B+G.$$
\end{prop}

\begin{defn}[ACSS, {cf. \cite[Definitions 5.4.2, 7.2.2, 7.2.3]{CHLX23}}]\label{defn: ACSS f-triple}
Let $(X,\Ff,B)$ be an lc foliated triple, $G\geq 0$ a reduced divisor on $X$, and $f: X\rightarrow Z$ a contraction. We say that $(X,\Ff,B;G)/Z$ is \emph{ACSS} if the following conditions hold:
\begin{enumerate}    
\item $(X,\Ff,B;G)/Z$ satisfies Property $(*)$.
\item $f$ is equidimensional.
\item There exists an $\Rr$-Cartier $\Rr$-divisor $D\geq 0$ on $X$, such that  $\Supp\{B\}\subset\Supp D$, and for any reduced divisor $\Sigma\geq f(G)$ such that $(Z,\Sigma)$ is log smooth, $$(X,B+D+G+f^*(\Sigma-f(G)))$$ 
      is qdlt (cf. \cite[Definition 35]{dFKX17}).
\item For any lc center of $(X,\Ff,B)$ with generic point $\eta$, over a neighborhood of $\eta$,
    \begin{enumerate}
      \item $\eta$ is the generic point of an lc center of $(X,\Ff,\lfloor B\rfloor)$, and
       \item $f: (X,B+G)\rightarrow (Z,f(G))$ is a toroidal morphism,
    \end{enumerate}
\end{enumerate}
 If $(X,\Ff,B;G)/Z$ is ACSS, then we say that $(X,\Ff,B)/Z$ and $(X,\Ff,B)$ are ACSS.
\end{defn}

\begin{defn}\label{def: relative ample model}
    Let $X',X,Z$ be normal quasi-projective varieties and $h: X'\rightarrow X$, $f: X'\rightarrow Z$ contractions. The \emph{core model} of $(h,f)$ associated with $(\bar h,\bar f)$ is the unique normal quasi-projective variety $\bar X$ up to isomorphism with two contractions $\bar h: \bar X\rightarrow X$ and $\bar f: \bar X\rightarrow Z$ satisfying the following.
    \begin{enumerate}
        \item  For any ample $\Rr$-divisor $A$ on $X$, $\bar h^*A$ is ample$/Z$.
        \item  There exists a contraction $g: X'\rightarrow \bar X$ such that $\bar h\circ g=h$ and $\bar f\circ g=f$.
    \end{enumerate}
The variety $\bar X$ is called the \emph{core model} of $(h,f)$ associated with $(\bar h,\bar f)$. Existence of core model is guaranteed by \cite[Definition-Lemma 3.1]{LMX24}.
\end{defn}

\begin{defn}[{\cite[Subsection 3.2]{LMX24}}]\label{deflem: simple model}
    Let $(X,\Ff,B)$ and $(X',\Ff',B')$ be two algebraically integrable foliated triples and $h: X'\rightarrow X$ a birational morphism. Let $f: X'\rightarrow Z$ be a contraction and $G$ a reduced divisor on $X'$. We say that $h: (X',\Ff',B';G)/Z\rightarrow (X,\Ff,B)$ is a \emph{simple modification} if the following conditions hold.
    \begin{itemize}
        \item $\Ff':=h^{-1}\Ff$ and $B':=h^{-1}_*(B\wedge\Supp B)+\Exc(h)^{\ninv}$, where $\Exc(h)^{\ninv}$ is the non-$\Ff'$-invariant part of $\Exc(h)$.
        \item $(X',\Ff',B')$ is lc.
        \item  $a(E,\Ff,B)\leq-\epsilon_{\Ff}(E)$ for any $h$-exceptional prime divisor $E$. 
        \item $K_{\Ff'}+B'\sim_ZK_{X'}+B'+G.$
        \item $(X',\Ff',B';G)/Z$ satisfies Property $(*)$.
    \end{itemize}
    We say that $h: (X',\Ff',B',G)/Z\rightarrow (X,\Ff,B)$ is 
    \begin{enumerate}
        \item an \emph{ACSS modification} if it is a simple modification and  $(X',\Ff',B';G)/Z$ is ACSS, and we say that $(X',\Ff',B';G)/Z$ is an \emph{ACSS model} of $(X,\Ff,B)$,
        \item a \emph{core modification} if it is a simple modification and $h^*A$ is ample$/Z$ for any ample $\Rr$-divisor $A$ on $X$, and we say that $(X',\Ff',B';G)/Z$ is an \emph{core model} of $(X,\Ff,B)$, and
        \item \emph{$\Qq$-factorial} if $X'$ is $\Qq$-factorial.
    \end{enumerate}
\end{defn}

\begin{lem}\label{lem: ACSS exist core mod equidim}
    Let $(X, \Ff,B)$ be an lc algebraically integrable triple. Then 
    \begin{enumerate}
        \item $(X,\Ff,B)$ has a $\Qq$-factorial ACSS model.
        \item There exists a core modification $\bar h: (\bar X,\bar \Ff, \bar B;\bar G)/Z\to (X, \Ff, B)$ associated with $\bar f: \bar X\to Z$ such that $\bar f$ is equidimensional.
    \end{enumerate}
\end{lem}

\begin{proof}
    (1) is \cite[Theorem 2.5.1]{CHLX23} or \cite[Theorem 3.10]{ACSS21}. We assume that $h: (X',\Ff',B';G)/Z \to (X,\Ff,B)$, associated with $f: X' \to Z$, is a $\Qq$-factorial ACSS model of $(X,\Ff,B)$. By \cite[Lemma 3.5]{LMX24}, the core model of $(h,f)$ associated with $(\bar h,\bar f)$ induces a core modification $(\bar X,\bar\Ff,\bar B;\bar G)/Z\rightarrow (X,\Ff,B)$. Since $f: X' \to Z$ is equidimensional, then $\bar f: \bar X \to Z$ is also equidimensional. (2) follows. 
\end{proof}

\begin{lem}\label{lem: core model is product}
Let $f: X\rightarrow Z$ be a contraction between normal quasi-projective varieties, and let $h_Z: Z'\rightarrow Z$ be a birational morphism. Let $h: X'\rightarrow X$ be a birational morphism and $f': X'\rightarrow Z'$ a contraction such that $f\circ h=h_Z\circ f'$. Let $\bar X$ be the normalization of the main component of $X\times_ZZ'$. Then $\bar X$ is the core model of $(h,f')$ associated with natural projections $(\bar h:\bar X\rightarrow X,\bar f: \bar X\rightarrow Z')$. In particular, for any ample$/Z$ $\Rr$-divisor $A$ on $X$, $h^*A$ is ample$/Z'$.
\end{lem}
\begin{proof}
By Stein factorization and Zariski's Main theorem, $\bar{h}$ and $\bar{f}$ are contractions. By the universal property of $X\times_Z Z'$, there exists a contraction $g: X'\rightarrow \bar X$ such that $\bar h\circ g=h$ and $\bar f\circ g=f$. For any ample$/Z$ $\Rr$-divisor $A$ on $X$, $\bar{h}^*A$ is ample over $Z'$. Hence by the uniqueness of core model, $\bar{X}$ is the core model.
\end{proof}

\subsection{Minimal model program}

For the definitions and notations of different types of minimal models for algebraically integrable foliations, we refer the reader to \cite[Definition 3.4]{LMX24}.

\begin{thm}\label{thm: ai foliation fiber gmm imply gmm}
Let $(X,\Ff,B)$ be a $\Qq$-factorial lc foliated triple such that $(X,\Ff,B;G)/Z$ satisfies Property $(*)$ for some $G\geq 0$ and contraction $f: X\rightarrow Z$. Let $F$ be a general fiber of $f$ and let $K_F+B_F:=(K_{\Ff}+B)|_F$. Assume that $(F,B_F)$ has a good minimal model. Then for any ample$/Z$ $\Rr$-divisor $A$ on $X$, we may run a $(K_{\Ff}+B)$-MMP$/Z$ with scaling of $A$ which terminates with a good minimal model of $(X,\Ff,B)/Z$.
\end{thm}
\begin{proof}
    Let $\Sigma_Z:=f(G)$ and let $0<\epsilon\ll 1$ be a real number. Let $\Delta:=B+G-\epsilon f^*\Sigma_Z$. Then any lc center of $(X,\Delta)$ is an lc center of $(X,B)$. By the definition of Property $(*)$, any lc center of $(X,\Delta)$ dominates $Z$. Moreover, since $(F,B_F)$ has a good minimal model and $K_F+B_F=(K_X+\Delta)|_F$, $(X,\Delta)$ is log abundant over $Z$. By \cite[Theorem 1.1]{Has23}, $(X,\Delta)$ has a good minimal model over $Z$.
    
    Possibly replacing $A$ with a general member in $|A|_{\mathbb R/Z}$, we may assume that $A\geq 0$ and $(X,\Delta+A)$ is lc. By \cite[Theorem 1.5]{HH20}, $(X,\Delta+\mu A)/Z$ has a good minimal model for any $\mu\in [0,1]$. By \cite[Lemma 5.4]{LMX24}, we may run a $(K_X+\Delta)$-MMP$/Z$ with scaling of $A$ which terminates with a good minimal model of $(X,\Delta)/Z$. By Proposition \ref{prop: weak cbf gfq}, this is also a $(K_{\Ff}+B)$-MMP$/Z$ with scaling of $A$ and terminates with a good minimal model of $(X,\Ff,B)/Z$.
\end{proof}

\subsection{Finite covers}

\begin{prop}[{\cite[Proposition 2.2]{CS20}, \cite[Lemmas 3.4, 4.3]{Dru21}, \cite[Proposition 3.7]{Spi20}}]\label{prop: cover formula}
    Let $(X,\Ff,B)$ be a foliated triple and let $h: X'\rightarrow X$ be a finite cover between normal varieties. Let $\Ff':=h^{-1}\Ff$. For any prime divisor $D$ on $X$, let $r_D$ be the ramification index of $h$ along $D$. Then:
    \begin{enumerate}
        \item Over the smooth locus of $X$, we have
          $$K_{\Ff'}:=h^*K_{\Ff}+\sum_D (r_D-1)D$$  
    where the sum runs through all non-$\Ff'$-invariant prime divisors on $X'$.
    \item Suppose that any codimension one component of the branch locus of $h$ is $\Ff$-invariant. Let $K_{\Ff'}+B':=h^*(K_{\Ff}+B)$, then:
    \begin{enumerate}
        \item If any component of $\Supp B$ is not $\Ff$-invariant, then the set of coefficients of $B'$ equals to the set of coefficients of $B$, and $h_*B'=(\deg h)B$.
        \item If $(X,\Ff,B)$ is lc (resp. klt, canonical, terminal), then $(X',\Ff',B')$ is also lc (resp. klt, canonical, terminal).
        \item If $(X',\Ff',B')$ is lc (resp. klt), then $(X,\Ff,B)$ is lc (resp. klt).
    \end{enumerate}
    \end{enumerate}
\end{prop}
\begin{proof}
(1) is \cite[Proposition 2.2]{CS20} or \cite[Lemma 3.4]{Dru21}; see \cite[Proposition 3.7]{Spi20} for the co-rank one case. (2.a) follows from (1) and that $h$ is unramified over the generic point of any component of $B$. (2.b) and (2.c) follow from the same lines of the proof of \cite[Lemma 4.3]{Dru21}. For the reader's convenience, we provide a full proof here.

Let $g: Y\rightarrow X$ be a birational morphism and $E$ a prime divisor on $Y$. Let $Y'$ be the normalization of the main component of $Y\times_{X}X'$ associated with $g': Y'\rightarrow X'$ and $h': Y'\rightarrow Y$. Let $E'$ be any prime divisor on $Y'$ that dominates $E$. Such prime divisor exists as $h'$ is finite. Let $\Ff_Y:=g^{-1}\Ff$ and $\Ff_{Y'}:=g'^{-1}\Ff'$.

Let $r_{E'}$ be the ramification index of $h'$ along $E'$. By (1), near the generic point of $E'$, we have
\begin{align*}
  K_{\Ff_{Y'}}&=g'^*(K_{\Ff'}+B')+a(E',\Ff',B')E'\\
  &=g'^*h^*(K_{\Ff}+B)+a(E',\Ff',B')E'\\
  &=h'^*g^*(K_{\Ff}+B)+a(E',\Ff',B')E'\\
  &=h'^*(K_{\Ff_Y}-a(E,\Ff,B)E)+a(E',\Ff',B')E'\\
  &=K_{\Ff_{Y'}}-\epsilon_{\Ff'}(E')(r_{E'}-1)E'-a(E,\Ff,B)r_{E'}E'+a(E',\Ff',B')E'.
\end{align*}
So
$$a(E',\Ff',B')+\epsilon_{\Ff'}(E')=r_{E'}(a(E,\Ff,B)+\epsilon_{\Ff'}(E'))=r_{E'}(a(E,\Ff,B)+\epsilon_{\Ff}(E)).$$
Therefore, $a(E,\Ff,B)\geq-\epsilon_{\Ff}(E)$ (resp. $>-\epsilon_{\Ff}(E)$) if and only if $a(E',\Ff',B')\geq-\epsilon_{\Ff'}(E')$ (resp. $>-\epsilon_{\Ff'}(E')$), and if $a(E,\Ff,B)\geq 0$ (resp. $>0$), then $a(E,\Ff,B)\geq 0$ (resp. $>0$). Moreover, $E$ is exceptional$/X$ if and only if $E'$ is exceptional$/X'$. By Lemma \ref{lem: equivalent definition lc and klt}, we get (2.c). 

By \cite[Lemma 2.22]{Kol13}, for any prime divisor $E'$ over $X'$, we may find $Y',Y,g,g',h'$ and $E$ as above. (2.b) follows.
\end{proof}

\begin{prop}[{cf. \cite[Proposition 3.5]{FS22}}]\label{prop: special cover}
        Let $(X,\Ff,B)$ be an algebraically integrable foliated triple. Assume that $(X,\Ff,B;G)/Z$ satisfies Property $(*)$ associated with an equidimensional contraction $f: X\rightarrow Z$. Then there exists a finite cover $h_Z: Z'\rightarrow Z$ between normal varieties satisfying the following. 
        
        Let $X'$ be the normalization of $X\times_ZZ'$, $h: X'\rightarrow X$ and $f': X'\rightarrow Z'$ the induced morphisms, $\Ff':=h^{-1}\Ff$, and
        $$K_{\Ff'}+B':=h^*(K_{\Ff}+B).$$
        Then:
        \begin{enumerate}
            \item   $f': (X',B')\rightarrow Z'$ is a locally stable family, and the set of coefficients of $B'$ equals to the set of coefficients of $B$.
            \item If $K_{\Ff}+B$ is ample$/Z$, then $f': (X',B')\rightarrow Z'$ is a stable family.
        \end{enumerate}
\end{prop}
\begin{proof}
(1) For any finite cover $h_Z: Z'\rightarrow Z$, the induced morphism $h: X'\rightarrow X$ is finite and any codimension one component of the branch locus of $h$ is not horizontal$/Z$, hence $\Ff$-invariant. (1) follows from \cite[Proposition 3.5]{FS22} and  Proposition \ref{prop: cover formula}(2.a). We note that although \cite[Proposition 3.5]{FS22} assumes $B\in \Qq$, the same line of proof also work for the case when $B\in\Rr$. Since $K_\Ff = K_X$ over the generic point of $Z$ and $B'$ is horizontal/$Z'$, $B'$ coincides with the boundary divisor defined in \cite[Proposition 3.5]{FS22}.


(2) Since $K_{\Ff}+B$ is ample$/Z$, $K_{\Ff'}+B'$ is ample$/Z$, so  $K_{\Ff'}+B'=K_{X'/Z'}+B'$ is ample$/Z'$. (2) follows from (1).
\end{proof}

\subsection{ACC for foliated pseudo-effective thresholds}

\begin{defn}
    Let $(X,\Ff,B)/U$ be a foliated triple and let $D\geq 0$ be an $\Rr$-Cartier $\Rr$-divisor on $X$. We define
    $$\pet(X/U,\Ff,B;D):=\inf\{t\geq 0,+\infty\mid K_{\Ff}+B+tD\text{ is pseudo-effective}/U, (X,\Ff,B+tD)\text{ is lc}\}$$
    to be the \emph{pseudo-effective threshold} of $D$ with respect to $(X,\Ff,B)/U$. For any positive integer $r$ and sets $\Ii\in [0,1]$, $\Ii'\in (0,+\infty)$, we define
\begin{align*}
\PET_{a.i.\fol}:=\Bigg\{\pet(X/U,\Ff,B;D)\Biggm|
    \begin{array}{r@{}l}
        \Ff\text{ is algebraically integrable},\\
        \rk\Ff=r,B\in\Ii,D\in\Ii'
    \end{array}\Bigg\}.
    \end{align*}
\end{defn}

\begin{thm}\label{thm: acc pet aif}
Let $r$ be a positive integer and $\Ii\subset [0,1]$, $\Ii'\subset (0,+\infty)$ two DCC sets. Then
$$\PET_{a.i.\fol}(r,\Ii,\Ii')$$
satisfies the ACC.
\end{thm}
\begin{proof}
We may assume that $1\in\Ii$.

Suppose that the theorem does not hold, then there exists a sequence of algebraically integrable foliated triples $\{(X_i,\Ff_i,B_i)/U_i\}_{i=1}^{+\infty}$ of rank $r$, $\Rr$-divisors $\{D_i\}_{i=1}^{+\infty}$, and a sequence of strictly increasing real numbers $\{t_i\}_{i=1}^{+\infty}$, such that $B_i,D_i\in\Ii$ and $t_i=\pet(X_i/U_i,\Ff_i,B_i;D_i)$. Let $t_0:=\lim_{i\rightarrow+\infty}t_i$. By the ACC for foliated lc thresholds for algebraically integrable foliations \cite[Theorem 1.1]{DLM23}, possibly passing to a subsequence, we may assume that $(X_i,\Ff_i,B_i+t_iD_i)$ is lc for each $i$ and $K_{\Ff_i}+B_i+(t_i-\epsilon)D_i$ is not pseudo-effective$/U$ for each $i$ and any positive real number $\epsilon$.

By Lemma \ref{lem: ACSS exist core mod equidim}, possibly replacing $(X_i,\Ff_i,B_i)$ with a $\Qq$-factorial ACSS model and replacing $D_i$ with its pullback, we may assume that at $(X_i,\Ff_i,B_i)$ is $\Qq$-factorial ACSS. In particular, $X_i$ is $\Qq$-factorial klt. By \cite[Theorem 1.4]{LMX24}, we may run a $(K_{\Ff_i}+B_i+(1-\frac{1}{i})t_iD_i)$-MMP$/U_i$ for each $i$ which terminates with a Mori fiber space $\phi_i: (X_i',\Ff_i',B_i'+(1-\frac{1}{i})t_iD_i')\rightarrow T_i$ of $(X_i,\Ff_i,B_i+(1-\frac{1}{i})t_iD_i)/U_i$, where $B_i',D_i'$ are the images of $B_i,D_i$ on $X_i'$ respectively. By \cite[Theorem 1.1]{DLM23}, possibly passing to a subsequence, we may assume that  $(X_i',\Ff_i',B_i'+t_0D_i')$ is lc for each $i$. 

Since $K_{\Ff_i}+B_i+t_iD_i$ is pseudo-effective$/U_i$, $K_{\Ff_i'}+B_i'+t_iD_i'$ is pseudo-effective$/U_i$. Since $\rho(X_i'/T_i)=1$, $K_{\Ff_i'}+B_i'+t_iD_i'$ is nef$/T_i$. Since $K_{\Ff_i'}+B_i'+(1-\frac{1}{i})t_iD_i'$ is anti-ample$/T_i$, $D_i'$ is ample$/T_i$, and there exist real number $\eta_i\in (1-\frac{1}{i},1]$ such that $K_{\Ff_i'}+B_i'+\eta_it_iD_i'\equiv_{T_i} 0$. Let $F_i$ be a general fiber of $\phi_i$. By the cone theorem for algebraically integrable foliations (cf. \cite[Theorem 3.9]{ACSS21}, \cite[Theorem 2.3.1]{CHLX23}), $F_i$ is tangent to $\Ff_i'$, so $K_{\Ff_i'}=K_{X_i'}$ over the generic fiber of $\phi_i$. In particular, $K_{\Ff_i'}|_{F_i}=K_{F_i}$.

Let $B_{F_i}:=B_i'|_{F_i}$ and $D_{F_i}:=D_i'|_{F_i}$, then $B_{F_i}\in\Ii,D_{F_i}\in\Ii'$, and
$$K_{F_i}+B_{F_i}+\eta_it_iD_{F_i}\equiv 0.$$
Since $(X_i',\Ff_i',B_i'+t_0D_i')$ is lc, $(F_i,B_{F_i}+\eta_it_iD_{F_i})$ is lc. Since $D_i'$ is ample$/T_i$, $D_{F_i}\not=0$. This contradicts the global ACC \cite[Theorem 1.5]{HMX14}.
\end{proof}

\begin{cor}\label{cor: 1 gap pet}
    Let $d$ be a positive integer and $\Ii$ a DCC set of real numbers. Then there exists a positive real number $\epsilon\in (0,1)$ depending only on $d$ and $\Ii$ satisfying the following. If $(X,\Ff,B)/U$ is an lc algebraically integrable foliated triple such that $\rk\Ff\leq d$, $B\in\Ii$, and $K_{\Ff}+B$ is big$/U$, then
    $$K_{\Ff}+B\wedge (1-\epsilon)\Supp B$$
    is big$/U$.
\end{cor}
\begin{proof}
We may assume that $1\in\Ii$. Let $h: (X',\Ff',B';G)/Z\rightarrow (X,\Ff,B)$ be a $\Qq$-factorial ACSS modification.

For any positive real number $\delta$, let $B_{\delta}'$ be the unique $\Rr$-divisor on $X'$ such that $\Supp B_{\delta}'=\Supp B'$, and for any component of $B'$, $\mult_DB_{\delta}'=1$ if $\mult_DB'\geq 1-\delta$, and  $\mult_DB_{\delta}'=\mult_DB'$ if $\mult_DB'<1-\delta$. By \cite[Theorem 1.1]{DLM23}, there exists a positive real number $\delta_0$, such that for any $\delta \in (0,\delta_0)$, $(X',\Ff',B_{\delta}')$ is lc. Since $K_{\Ff}+B$ is big$/U$, $K_{\Ff'}+B'$ is big$/U$, so $K_{\Ff'}+B_{\delta}'$ is big$/U$ for any $\delta\in (0,1)$.

By Theorem \ref{thm: acc pet aif}, there exists a real number $\epsilon\in (0,\delta_0)$ depending only on $d$ and $\Ii$, such that $K_{\Ff'}+\{B_{\delta}'\}+(1-\epsilon)\lfloor B_{\delta}'\rfloor$ is big$/U$ for any $\delta\in (0,\delta_0)$. In particular,
$$K_{\Ff'}+\{B_{\epsilon}'\}+(1-\epsilon)\lfloor B_{\epsilon}'\rfloor=K_{\Ff'}+B'\wedge (1-\epsilon)\Supp B'$$
is big$/U$. Thus
$$K_{\Ff}+B\wedge (1-\epsilon)\Supp B=h_*(K_{\Ff'}+B'\wedge (1-\epsilon)\Supp B')$$
is big$/U$.
\end{proof}

\section{Proof of the main theorems}

 We first prove the following special case of Theorem \ref{thm: main theorem foliation with boundary} under the additional assumption that $\Ff$ is induced by a contraction and $(X,\Ff,B)$ is already relatively stable.

\begin{thm}\label{thm: main theorem foliation already induced case}
Let $d$ be a positive integer, $\Ii\subset [0,1]$ a DCC set, and $v$ a positive real number. Then there exists a positive real number $\delta$ and a discrete set $\Ii_0$ depending only on $d,\Ii$ and $v$ satisfying the following. Assume that 
\begin{enumerate}
\item $(X,\Ff,B)$ is a projective lc foliated triple,
\item $\Ff$ is induced by a contraction $f: X\rightarrow Z$,
\item $(X_g,B_g)$ is a general fiber of $f: (X,B)\rightarrow Z$,
\item $\dim X_g=d$ and the coefficients of $B$ belong to $\Ii$,
\item $K_{\Ff}+B$ is ample$/Z$ and $\vol(K_{X_g}+B_g)=v$.
\end{enumerate}
Then $\vol(K_{\Ff}+B)\in\Ii_0$.
\end{thm}

\begin{proof} 
We may assume that $\vol(K_{\Ff}+B)>0$. By \cite[6.8.4]{Kol23}, we may assume that $\Ii$ is a finite set. Possibly replacing $\Ii$ with $\Ii_+$, we may assume that $\Ii=\Ii_+$. By Theorem \ref{kp17 6.19}, there exists a reduced scheme $\Ss$ and a family $\ff: (\Xx,\Bb)\rightarrow\Ss$ such that $\ff\in\mathfrak{S}(d,\Ii,v)(\Ss)$ and the induced moduli map $\Ss\rightarrow M(d,\Ii,v)$ is finite with degree $s$. Note that $\dim\Ss$ and $s$ depend only on $d,\Ii$ and $v$.

\medskip

\noindent\textbf{Step 1}. In this step, we reduce to the case when $Z$ is smooth, $(X,\Ff,B;G)/Z$ satisfies Property $(*)$ for some $G\geq 0$, $X\to Z$ is equidimensional, and there exists a morphism $Z\rightarrow\Ss$ and a non-empty open subset $Z^0\subset Z$, such that the induced morphism $(X,B)\times_ZZ^0\rightarrow Z^0$ is a family of $(d,\Ii,v)$-stable pairs with induced moduli map $Z^0\rightarrow\Ss\rightarrow M(d,\Ii,v)$.

 By our assumption, $(X_g,B_g)$ is a $(d,\Ii,v)$-stable pair. Therefore, there exists a non-empty open subset $Z^0\subset Z$ such that $f^0: (X^0,B^0)\rightarrow Z^0$ is a family of $(d,\Ii,v)$-stable pairs, where $X^0:=X|_{f^{-1}(Z^0)}$, $B^0:=B|_{X^0}$, and $f^0:=f|_{X^0}$. Thus there exists a moduli map $\phi_0: Z^0\rightarrow M(d,\Ii,v)$ which naturally extends to a rational map $\phi: Z\dashrightarrow M(d,\Ii,v)$. 

Since $X$ is normal and $f$ is a contraction, $Z$ is normal. Let $h_Z: Z'\rightarrow Z$ be a resolution such that the induced map $\phi':=\phi\circ h_Z: Z'\dashrightarrow M(d,\Ii,v)$ is a morphism. Let $X'$ be the normalization of the main component of $X\times_ZZ'$ with induced morphisms $h: X'\rightarrow X$ and $f': X'\rightarrow Z'$. Let $\Ff':=h^{-1}\Ff$ and 
$$K_{\Ff'}+B':=h^*(K_{\Ff}+B).$$
Since $(\Ff,B)$ is lc, $B'\in\Ii$. Since $h$ is birational, we have
$$\vol(K_{\Ff'}+B')=\vol(K_{\Ff}+B).$$


Let $\bar Z$ be the normalization of $Z'\times_{M(d,\Ii,v)}\Ss$, $h_Z': \bar Z\rightarrow Z'$ the induced morphism, and $\bar X$ the normalization of the main component of $X'\times_{Z'}\bar Z$ with induced morphisms $h': \bar X\rightarrow X'$ and $\bar f: \bar X\rightarrow\bar Z$. Since $h_Z'$ is finite, $h'$ is finite. Let $\bar\Ff:=h'^{-1}\Ff'$ and let
$$K_{\bar \Ff}+\bar B=h'^*(K_{\Ff'}+B').$$
We note that $\bar f$ is induced by $\bar X\to \bar Z$. Since the branch locus of $h_Z'$ cannot contain the whole $\bar{Z}$ and unramification is stable under change base, the branch locus of $h'$ is vertical over $\bar{Z}$, thus does not contain any non-$\Ff'$-invariant divisor. By Proposition \ref{prop: cover formula}, $(\bar X,\bar\Ff,\bar B)$ is lc. Let $(\bar X_g,\bar B_g)$ be a general fiber of $\bar f: (\bar X,\bar B)\rightarrow \bar{Z}$. By our construction, $\dim\bar X_g=d$, $\bar B\in\Ii$, $K_{\bar\Ff}+\bar B$ is ample$/\bar Z$, and $\vol(K_{\bar X_g}+\bar B_g)=v$, and there exists a non-empty open subset $\bar Z^0$ of $\bar Z$ such that the induced morphism $(\bar X,\bar B)\times_{\bar Z}\bar Z^0\rightarrow \bar Z^0$ is a family of $(d,\Ii,v)$-stable pairs with induced moduli map $\bar Z^0\rightarrow \Ss\rightarrow M(d,\Ii,v)$. We have
$$\vol(K_{\bar\Ff}+\bar B)=\deg h'\vol(K_{\Ff'}+B')=s\vol(K_{\Ff'}+B').$$
Therefore,
$$\vol(K_{\Ff}+B)=\frac{1}{s}\vol(K_{\bar\Ff}+\bar B).$$

Now we let $(\bar X',\bar\Ff',\bar B';\bar G')/\bar Z'\rightarrow (\bar X,\bar\Ff,\bar B)$ be a core model of $(\bar X,\bar\Ff,\bar B)$ whose existence is guaranteed by Lemma \ref{lem: ACSS exist core mod equidim}. Then $\bar Z'$ is smooth, $\bar X' \to \bar Z'$ is equidimensional, and $\bar X'$ is the normalization of the main component of $\bar X\times_{\bar Z}\bar Z'$. Let $\bar Z^{0'}$ be the inverse image of the intersection of $\bar Z^0$. We may replace $(X,\Ff,B)\rightarrow Z^0$ with $(\bar X',\bar\Ff',\bar B')\rightarrow\bar Z'$ and $Z$ with $\bar Z'$, and conclude \textbf{Step 1}. 

\medskip

\noindent\textbf{Step 2}. Let $\ff_Z: (\Xx_Z,\Bb_Z)\rightarrow Z$ be the base change of $\ff: (\Xx,\Bb)\rightarrow\Ss$.  In this step we 
reduce to the case when $f: (X,B)\rightarrow Z$ is $\ff_Z: (\Xx_Z,\Bb_Z)\rightarrow Z$. For simplicity, from now on we forget all new notations introduced in \textbf{Step 1} except $Z^0$.

By \textbf{Step 1}, there exists a non-empty open subset $Z^0\subset Z$ such that $f^0:(X^0,B^0):=(X,B)\times_ZZ^0\rightarrow Z^0$ is a family of $(d,\Ii,v)$-stable pairs and the moduli map $Z^0\rightarrow M(d,\Ii,v)$ factors through $\Ss$. Since the generic fiber of $f: (X,B)\rightarrow Z$ is normal, the generic fiber of $f^0:(X^0,B^0)\rightarrow Z^0$ is normal. By Proposition \ref{prop: isom functor is finite}, the generic fiber of $\ff_Z: (\Xx_Z,\Bb_Z)\rightarrow Z$ is normal. We note that $\Xx_Z$ is normal as $Z$ is smooth and the generic fiber of $\ff_Z$ is normal (cf. \cite[10.78]{Kol23}). 

By Proposition \ref{prop: special cover}, there exists a finite cover $h_Z: Z'\rightarrow Z$ bewteen normal varieties satisfying the following: Let $X'$ be the normalization of $X\times_ZZ'$, $h: X'\rightarrow X$ and $f': X'\rightarrow Z'$ the induced morphisms, and let $K_{\Ff'}+B':=h^*(K_{\Ff}+B)$, then $f': (X',B')\rightarrow Z'$ is a stable family of $(d,\Ii,v)$-stable pairs.

 Let $\ff': (\Xx',\Bb')\rightarrow Z'$ be the base change of $\ff: (\Xx,\Bb)\rightarrow\Ss$, then it is also the base change of $\ff_Z: (\Xx_Z,\Bb_Z)\rightarrow Z$. By our assumption, over a non-empty open subset $Z'^0$ of $Z'$, the moduli map $Z'^0\rightarrow M(d,\Ii,v)$ for the family $(X',B')\times_{Z'}Z'^0\rightarrow Z'^0$ factors through $\Ss$. By the separatedness of $M(d,\Ii,v)$, the moduli map $Z'\rightarrow M(d,\Ii,v)$ of the family $f': (X',B')\rightarrow Z'$ is the same as the moduli map $Z'\rightarrow M(d,\Ii,v)$ of the family $\ff': (\Xx',\Bb')\rightarrow Z'$. 
 
 Let $\mathfrak{h}: \Xx'\rightarrow\Xx_Z$ be the induced morphism. Then
$\deg\mathfrak{h}=\deg h_Z=\deg h$. By Lemma \ref{lem: over same base lsf volume are same},
\begin{align*}
\deg\mathfrak{h}\cdot\vol(K_{\Xx_Z/Z}+\Bb)=&\vol(K_{\Xx'/Z'}+\Bb')=\vol(K_{X'/Z'}+B')\\
=&\vol(K_{\Ff'}+B')=\deg h\cdot\vol(K_{\Ff}+B).
\end{align*}
Thus
$$\vol(K_{\Ff}+B)=\vol(K_{\Xx_Z/Z}+\Bb_Z).$$
Since the generic fiber of $\ff_Z: (\Xx_Z,\Bb_Z)\rightarrow Z$ is normal, conditions (1-5) hold for $\ff_Z: (\Xx_Z,\Bb_Z)\rightarrow Z$ and the foliation induced by $\ff_Z$. So we may replace $f: (X,B)\rightarrow Z$ with $\ff_Z: (\Xx_Z,\Bb_Z)\rightarrow Z$ and conclude \textbf{Step 2}.

\medskip

\noindent\textbf{Step 3}. In this step we conclude the proof. Let $T$ be the normalization of the image of $Z$ in $\Ss$ and let $f_T: (X_T,B_T)\rightarrow T$ be the base change of $\ff: (\Xx,\Bb)\rightarrow\Ss$. We note that $X_T$ is normal as $T$ is normal and the general fiber of $f$ is normal (cf. \cite[10.78]{Kol23}).

Since $K_{X/Z}+B=K_{\Ff}+B$ is big, by Lemma \ref{lem: non maximal variation not big}, $f: (X,B)\rightarrow Z$ is of maximal variation. Thus the moduli map $\phi: Z\rightarrow M(d,\Ii,v)$ is a generically finite morphism. Therefore, the induced morphism $\tau: Z\rightarrow T$ is generally finite. We let $l:=\deg\tau$, then
$$l\vol(K_{X_T/T}+B_T)=\vol(K_{X/Z}+B).$$


Since ampleness is an open condition, by \cite[Theorem 5.6]{HLS19}, there exist positive real numbers $a_1,\dots,a_k$ such that $\sum_{i=1}^ka_i=1$, $\Bb=\sum_{i=1}^k a_i\Bb_i$, such that $\ff: (\Xx,\Bb_i)\rightarrow\Ss$ is a stable family for each $i$, and each $\Bb_i$ is a $\Qq$-divisor. We let $f_T: (X_T,B_{i,T})\rightarrow T$ be the base change of $\ff: (\Xx,\Bb_i)\rightarrow\Ss$, then for each $i$, $f_T: (X_T,B_{i,T})\rightarrow T$ is a stable family. 
Let $I$ be a positive integer, depending only on $(\Xx,\Bb)\rightarrow\Ss$ hence depending only on $d,\Ii,v$, such that $I(K_{\Xx/\Ss}+\Bb_i)$ is Cartier for each $i$. Then $H_i:=I(K_{X_T/T}+B_{i,T})$ is Cartier for each $i$ because it is the pullback of $I(K_{\Xx/\Ss}+\Bb_i)$ to $X_T$.

Note that the generic fiber of $f_T$ is normal as $X,T$ are normal. By Proposition \ref{prop: maximal variation big and nef}, $H_i$ is big and nef for each $i$. Therefore,
$$\vol(K_{X_T/T}+B_T)=\frac{1}{I^{\dim X_T}}\left(\sum_{i=1}^ka_iH_i\right)^{\dim X_T}.$$
Since $\dim X_T\leq\dim\Xx=d+\dim\Ss$, $\vol(K_{X_T/T}+B_T)$ belongs to a discrete set. Thus
$$\vol(K_{\Ff}+B)=\vol(K_{X/Z}+B)=l\vol(K_{X_T/T}+B_T)$$
belongs to a discrete set and we conclude the proof.
\end{proof}

\begin{proof}[Proof of Theorem \ref{thm: main theorem foliation with boundary}]
We may assume that $1\in\Ii$ and $K_{\Ff}+B$ is big. Let $h: (X',\Ff',B';G)/Z\rightarrow (X,\Ff,B)$ be a $\Qq$-factorial ACSS modification of $(X,\Ff,B)$ associated with contraction $f: X'\rightarrow Z$. Let $L'$ be the leaf correspond to $L$, and $h_{L'}: L'\to L$ the induced morphism. Then $L'$ is a general fiber of $X'\to Z$ and $B'\in\Ii$.

Let $K_{L'}+B_{L'}:=(K_{\Ff'}+B')|_{L'}$, then
$$K_{L'}+B_{L'}=h_{L'}^*(K_L+B_L),$$
so $\vol(K_{L'}+B_{L'})=v$. Moreover, 
$$\vol(K_{\Ff'}+B')=\vol(K_{\Ff}+B)>0.$$

We first prove (2). Since $(X_L,B_L)$ has a good minimal model, then by \cite[Lemma 2.14]{Has19}, $(L',B_{L'})$ has a good minimal model. By Theorem \ref{thm: ai foliation fiber gmm imply gmm}, we may take the ample model $g:X'\dashrightarrow \bar X$ of $K_{\Ff'}+B'/Z$. Let $\bar\Ff:=g_*\Ff'$, and $\bar B:=g_*B'$. Since $K_{\Ff'}+B'$ is big, $g$ is birational. Thus $(\bar X,\bar\Ff,\bar B)$ is lc, and we have
$$\vol(K_{\bar\Ff}+\bar B)=\vol(K_{\Ff'}+B')=\vol(K_{\Ff}+B).$$

Let $\bar L$ be the images of $L'$ on $\bar X$, and $B_{\bar L}:=\bar B|_{\bar L}$. Then $g|_{L'}: L'\dashrightarrow \bar{L}$ is the ample model of $K_{L''}+B_{L''}$. We have $$\vol(K_{\bar L}+B_{\bar L})=\vol(K_{L'}+B_{L'})=v,$$
and $K_{\bar L}+B_{\bar L}$ is ample. By Theorem \ref{thm: main theorem foliation already induced case}, there exists a finite set $\Ii_0$ depending only on $d,\Ii$ and $v$, such that $\vol(K_{\Ff}+B)=\vol(K_{\bar\Ff}+\bar B)\in\Ii_0$. This implies (2). 



Now we prove (1). Let $\epsilon\in (0,1)$ be the number as in Corollary \ref{cor: 1 gap pet} which only depends on $d$ and $\Ii$. Let $\tilde B':=B'\wedge (1-\epsilon)\Supp B'$ and let $\tilde B_{L'}:=\tilde B'|_{L'}$. Since $X'$ is $\Qq$-factorial klt and $(X',B')$ is qdlt, $(X',\tilde B')$ is klt. Thus there are finitely many prime divisors over $X'$ whose discrepancies with respect to $(X',\tilde B')$ are $\leq -1+\epsilon$ and dominate $Z$. 

Let $p: Y\rightarrow X'$ be an extraction of all these divisors such that $Y$ is $\Qq$-factorial, whose existence is guaranteed by \cite[Corollary 1.4.3]{BCHM10}, $\Ff_Y:=p^{-1}\Ff'$, 
$$K_{\Ff_Y}+B_Y:=p^*(K_{\Ff'}+B'),$$
and
$$\tilde B_Y:=B_Y\wedge (1-\epsilon)\Supp B_Y.$$
Since $K_{\Ff'}=K_{X'}$ over the generic point of $Z$,
$$K_Y+B_Y=p^*(K_{X'}+B').$$
Let $L_Y$ be a general fiber of the induced contraction $Y\rightarrow Z$ and let $K_{L_Y}+B_{L_Y}:=(K_{\Ff_Y}+B_Y)|_{L_Y}$.
Since $K_{\Ff_Y}+B_Y$ is big, by our construction of $\epsilon$, $K_{\Ff_Y}+\tilde B_Y$ is big. Since $K_L+B_L$ is big, $K_{L_Y}+B_{L_Y}=p|_{L_Y}^*(K_L+B_L)$ is big. By the construction of $\epsilon$, $K_{L_Y}+\tilde B_{L_Y}$ is big, where 
$$\tilde B_{L_Y}:=\tilde B_Y|_{L_Y}.$$
By our construction, $(Y,\tilde B_Y)$ is $\epsilon$-lc over the generic point of $Z$ and $\tilde B_Y\in\Ii\cup\{1-\epsilon\}$. Thus $(L_Y,\tilde B_{L_Y})$ is $\epsilon$-lc and
$$0<\vol(K_{L_Y}+\tilde B_{L_Y})\leq\vol(K_{L_Y}+B_{L_Y})=\vol(K_L+B_L)=v.$$
By Theorem \ref{thm: acc pet aif}, there exists a positive integer $m$ depending only on $d$ and $\Ii$, such that $K_{\Ff_Y}+\frac{1}{m}\lfloor m\tilde B_Y\rfloor$ and $K_{L_Y}+\frac{1}{m}\lfloor m\tilde B_{L_Y}\rfloor$ are big. Then $\left(L_Y,\frac{1}{m}\lfloor m\tilde B_{L_Y}\rfloor\right)$ is $\epsilon$-lc and $$0<\vol\left(K_{L_Y}+\frac{1}{m}\lfloor m\tilde B_{L_Y}\rfloor\right)\leq \vol(K_{L_Y}+\tilde B_{L_Y})\leq v.$$
By \cite[Theorem 1.3]{Fil20}, $\vol\left(K_{L_Y}+\frac{1}{m}\lfloor m\tilde B_{L_Y}\rfloor\right)$ belongs to a finite set depending only on $d,\Ii$ and $v$. By \cite[Theorem 1.2]{BCHM10}, $\left(L_Y,\frac{1}{m}\lfloor m\tilde B_{L_Y}\rfloor\right)$ has a good minimal model. Since $K_{\Ff_Y}+\frac{1}{m}\lfloor m\tilde B_Y\rfloor$ is big, by (2), there exists a positive real number $\delta$ depending only on $d,\Ii$ and $v$, such that 
$$\vol(K_{\Ff}+B)=\vol(K_{\Ff_Y}+B_Y)\geq\vol(K_{\Ff_Y}+\tilde B_Y)\geq\vol\left(K_{\Ff_Y}+\frac{1}{m}\lfloor m\tilde B_Y\rfloor\right)\geq\delta.$$
The theorem follows.
\end{proof}

\begin{cor}\label{cor: main theorem foliation with boundary klt case}
Let $d$ be a positive integer, $v$ a positive real number, and $\Ii\subset [0,1]$ a DCC set. Then there exists a positive real number $\delta$ and a discrete set $\Ii_0$ depending only on $d,v$ and $\Ii$ satisfying the following. Assume that 
\begin{itemize}
    \item $(X,\Ff,B)$ is a projective klt algebraically integrable foliation,
    \item  $L$ is the normalization of a general leaf of $\Ff$ and $B_L:=B|_L$,
    \item $\dim L=d$ and the coefficients of $B$ belong to $\Ii$, and
    \item $\vol(K_L+B_L)=v$.
\end{itemize}
Then  $\vol(K_{\Ff}+B)\in\Ii_0$. In particular, of $K_{\Ff}+B$ is big, then $\vol(K_{\Ff}+B)\geq\delta$.
\end{cor}
\begin{proof}

    Let $h: (X',\Ff',B';G)/Z\rightarrow (X,\Ff,B)$ be a $\Qq$-factorial ACSS modification of $(X,\Ff,B)$. Since $(X,\Ff,B)$ is klt and $a(E,\Ff,B)=-\epsilon_{\Ff}(E)$ for any prime $h$-exceptional divisor $E$, $h$ does not extract any non-$\Ff'$-invariant divisor, and $(X',B')$ is klt. Let $L'$ be the leaf correspond to $L$, and $h_{L'}: L'\to L$ the induced morphism. Then $h_{L'}$ is small and hence an isomorphism. Thus $(L',B_{L'}:=B'|_L)$ is klt which implies that $(L,B_L)$ is klt. By \cite[Theorem 1.2]{BCHM10}, $(L,B_L)$ has a good minimal model. The corollary follows from Theorem \ref{thm: main theorem foliation with boundary}.
\end{proof}

\begin{proof}[Proof of Theorem \ref{thm: main theorem foliation no boundary}]
It is a special case of Corollary \ref{cor: main theorem foliation with boundary klt case}.
\end{proof}

\begin{proof}[Proof of Theorem \ref{thm: main theorem lsf with boundary}]
Since birational base change does not change the volume of the relative log canonical divisor nor the general log fiber, we may take a normalization of $Z$ and assume that $Z$ is normal. Now Theorem \ref{thm: main theorem lsf with boundary} is a special case of Theorem \ref{thm: main theorem foliation with boundary}. 
\end{proof}

\begin{proof}[Proof of Theorem \ref{thm: main theorem lsf no boundary}] It is a consequence of Theorem \ref{thm: main theorem lsf with boundary} and Lemma \ref{lem: non maximal variation not big}.
\end{proof}

\begin{rem}
Generalized pairs \cite{BZ16} and generalized foliated quadruples \cite{LLM23} play crucial roles in recent developments of birational geometry and foliation theory. We expect Theorem \ref{thm: main theorem lsf no boundary} and \ref{thm: main theorem foliation with boundary} to hold for locally stable families of (NQC) generalized pairs and generalized foliated quadruples respectively. This is because the minimal model programs for NQC generalized pairs and generalized foliated quadruples are known \cite{BZ16,HL22,HL23,LX23,Xie22}, and the volume for NQC generalized pairs is well-studied \cite{Bir21b}. However, since we do not have a robust moduli theory for lc generalized pairs \cite{BH22} and the moduli theory for klt generalized pairs is not yet established, we cannot prove these analogous results at this time.
\end{rem}

\end{document}